\documentclass{article}%
\usepackage{amsmath}
\usepackage{graphicx}
\usepackage{amsfonts}
\usepackage{amssymb}%
\newcommand{\ab}{\linebreak[0]} 
\newcommand{\sep}{\vspace{-3pt} \begin{center}
{\mathversion{normal}
$\infty \mspace{-5.5mu} \infty \mspace{-5.5mu}
\infty \mspace{-5.5mu} \infty \mspace{-5.5mu}
\infty \mspace{-5.5mu} \infty \mspace{-5.5mu}
\infty \mspace{-5.5mu} \infty$
}
\end{center} \vspace{-3pt}}
\begin{document}

\author{David J. Pengelley\thanks{Dedicated to the memory of my parents, Daphne and Ted Pengelley, who inspired 
a love of history.}
\\Department of Mathematical Sciences\\New Mexico State University\\Las Cruces, NM 88003, USA}
\title{Dances between continuous and discrete: \\Euler's summation formula\thanks{Based on a talk given at the Euler 2K+2
conference, Rumford, Maine, 2002.}}
\date{May 5, 2007}
\maketitle

\section{Introduction}

Leonhard Euler (1707--1783) discovered his powerful \textquotedblleft
summation formula\textquotedblright\ in the early 1730s. He used it in 1735 to
compute the first $20$ decimal places for the precise sum of all the
reciprocal squares --- a number mathematicians had competed to determine ever
since the surprising discovery that the alternating sum of reciprocal odd
numbers is $\pi/4$. This reciprocal squares challenge was called the
\textquotedblleft Basel problem,\textquotedblright and Euler achieved his
$20$-place approximation using only a few terms from his diverging
summation formula. In contrast, if sought as a simple partial sum of the
original slowly converging series, such accuracy would require 
more than $10^{20}$ terms. With his approximation, Euler probably
became convinced that the sum was $\pi^{2}/6$, which spurred his
first solution of the Basel problem in the same year
\cite[volume\ 16, section 2, pp.\ VIIff, volume\ 14]{euleroo}%
\cite{masters,weil}.

We are left in awe that just a few terms of a diverging formula can so closely
approximate this sum. Paradoxically, Euler's formula, even though it usually
diverges, provides breathtaking approximations for partial and
infinite sums of many slowly converging or diverging series. My goal here is
to explore Euler's own mature view of the summation formula and a few of his
more diverse applications, largely in his own words from the 
\emph{Institutiones Calculi Differentialis} \emph{(Foundations of Differential
Calculus) }of 1755. I hope that readers
will be equally impressed at some of his other applications.

In the \emph{Calculi Differentialis}, Euler connected his summation formula to
Bernoulli numbers and proved the sums of powers formulas that Jakob Bernoulli
had conjectured. He also applied the formula to harmonic partial sums and the
related gamma constant, and to sums of logarithms, thereby approximating large
factorials (Stirling's asymptotic approximation) and binomial coefficients
with ease. He even made an approximation of $\pi$ that he himself commented
was hard to believe so accurate for so little work. Euler was a wizard at
finding these connections, at demonstrating patterns by generalizable example,
at utilizing his summation formula only \textquotedblleft until it begins to
diverge,\textquotedblright and at determining the relevant \textquotedblleft
Euler-Maclaurin constant\textquotedblright\ in each application. 
His work also inaugurated study of the zeta function
\cite{ayoub,schup}. Euler's accomplishments throughout this entire arena are
discussed from different points of view in many modern books
\cite{dunhameuler}\cite[pp.\ 119--136]{goldstine}\cite[II.10]{hairer}%
\cite[chapter\ XIII]{hardy}\cite[p.\ 197ff]{hildebra}\cite[chapter XIV]{knopp}
\cite[p.\ 184, 257--285]{weil}\cite[p.\ 338ff]{young}.

Euler included all of these discoveries and others in beautifully unified form
in Part Two\footnote{Part One has recently appeared in English translation
\cite{eulerfdcblanton}, but not Part Two.} of the \emph{Calculi Differentialis
}\cite[volume 10]{euleroo}\cite{eulertransl}, portions of which I have
translated for an undergraduate course based on original sources
\cite{recoveringmotivation,411course,webresource}, and for selective
inclusion\footnote{See \cite{instcalcdiff} for my most extensive translation
from Euler's Part Two (albeit more lightly annotated).} in a companion book
built around annotated primary sources \cite{masters}. The chapter \emph{The
Bridge Between Continuous and Discrete} \cite{masters,webresource} follows the
entwining of the quest for formulas for sums of numerical powers with the
development of integration, via sources by Archimedes, Fermat,
Pascal, Jakob Bernoulli, and finally from Euler's \emph{Calculi
Differentialis}.\ I have also written an article \cite{bridge} providing an
independent exposition of this broader story.

Here I will first discuss the Basel problem and briefly outline the
progression of ideas and sources that led to the connection in Euler's
work between it and sums of powers. Then I will illustrate a
few of Euler's achievements with his summation formula via selected
translations. I present Euler's derivation of the formula, discuss his
analysis of the resulting Bernoulli numbers, show his application to sums of
reciprocal squares, to large factorials and binomial coefficients, and mention
other applications. A more detailed treatment can be found in \cite{masters}.
I will also raise and explore the question of whether large factorials can be
determined uniquely from Euler's formula.

\section{The Basel problem}

In the 1670s, James Gregory (1638--1675) and Gottfried Leibniz (1646--1716)
discovered that
\[
1-\frac{1}{3}+\frac{1}{5}-\frac{1}{7}+\cdots=\frac{\pi}{4},
\]
as essentially had the mathematicians of Kerala in southern India two
centuries before \cite[pp.\ 493ff,527]{katz}. Because, aside from geometric
series, very few infinite series then had a known sum, this remarkable result
enticed Leibniz and the Bernoulli brothers Jakob (1654--1705) and Johann
(1667--1748) to seek sums of other series, particularly the reciprocal
squares
\[
\frac{1}{1}+\frac{1}{4}+\frac{1}{9}+\frac{1}{16}+\frac{1}{25}+\cdots=\text{
?,}%
\]
a problem first raised by Pietro Mengoli (1626--1686) in 1650. Jakob expressed
his eventual frustration at its elusive nature in the comment ``\emph{If
someone should succeed in finding what till now withstood our efforts and
communicate it to us, we shall be much obliged to him}'' \cite[p.\ 345]{young}.

Euler proved that the sum is exactly $\pi^{2}/6$, in part by a broadening of
the context to produce his \textquotedblleft summation
formula\textquotedblright\ for $\sum_{i=1}^{n}f(i)$, with $n$ possibly
infinite. His new setting thus encompassed both the Basel problem, $\sum
_{i=1}^{\infty}1/i^{2}$, and the quest for closed formulas for sums of powers,
$\sum_{i=1}^{n}i^{k}\approx\int_{0}^{n}x^{k}\,dx$, which had been sought since
antiquity for area and volume investigations. The summation formula helped
Euler resolve both questions. This is a fine pedagogical illustration of how
generalization and abstraction can lead to the combined solution of seemingly
independent problems.

\section{Sums of powers and Euler's summation formula: historically
interlocked themes}

Our story (told more completely elsewhere \cite{masters,bridge}) begins in
ancient times with the Greek approximations used to obtain 
areas and volumes by the method of exhaustion. The Pythagoreans (sixth century
B.C.E.) knew that
\[
1+2+3+\cdots+n=\frac{n(n+1)}{2},
\]
and Archimedes (third century B.C.E.) proved an equivalent to our modern
formula
\[
1^{2}+2^{2}+3^{2}+\cdots+n^{2}=\frac{n(n+1)(2n+1)}{6},
\]
which he applied to deduce the area inside a spiral: \textquotedblleft%
\emph{The area bounded by the first turn of the spiral and the initial line is
equal to one-third of the first circle\textquotedblright} \cite[Spirals]{arch}.

Summing yet higher powers was key to computing other areas and volumes, and one
finds the formula for a sum of cubes in work of Nicomachus of Gerasa (first
century B.C.E.), \={A}ryabha\d{t}a in India (499 C.E.), and al-Karaj\={\i} in
the Arab world (circa 1000) \cite{boyersums}\cite[p.\ 68f]{heathmanual}%
\cite[p.\ 212f,251ff]{katz}. The first evidence of a general relationship
between various exponents is in the further Arabic work of Ab\={u} `Al\={\i}
al-\d{H}asan ibn al-Haytham (965--1039), who needed a formula for sums of
fourth powers to find the volume of a paraboloid of revolution. He discovered
a doubly recursive relationship between exponents \cite[p.\ 255f]%
{katz}.

By the mid-seventeenth century Pierre de Fermat (1601--1665) and Blaise Pascal
(1623--1662) had realized the general connection between the figurate
(equivalently binomial coefficient) numbers and sums of powers, motivated by
the drive to determine areas under \textquotedblleft higher
parabolas\textquotedblright\ (i.e., $y=x^{k}$) \cite[p.\ 481ff]{katz}. Fermat
called the sums of powers challenge \textquotedblleft\emph{what is perhaps the
most beautiful problem of all arithmetic,}\textquotedblright and he claimed a
recursive solution using figurate numbers. Pascal used binomial expansions and
telescoping sums to obtain the first simply recursive relationship between
sums of powers for varying exponents \cite{boyersums}.

Jakob Bernoulli, during his work in the nascent field of probability, was the
first to conjecture a general pattern in sums of powers formulas,
simultaneously introducing the Bernoulli numbers into mathematics\footnote{The
evidence suggests that around the same time, Takakazu Seki (1642?--1708) in
Japan also discovered the same numbers \cite{shen,yosida}.}. In his
posthumous book of 1713, \emph{The Art of Conjecturing} \cite[volume 3,
pp.\ 164--167]{bernoulli}, appears a section on \emph{A Theory of Permutations
and Combinations}. Here one finds him first list the formulas for \emph{Sums
of Powers} up to exponent ten (using the notation $\int$ for the discrete sum
from $1$ to $n$), and then claim a pattern, to wit\footnote{Bernoulli's
asterisks in the table indicate missing monomial terms. Also, there is an
error in the original published Latin table of sums of powers formulas. The
last coefficient in the formula for $\int n^{9}$ should be $-\frac{3}{20}$,
not $-\frac{1}{12}$; we have corrected this here.}:%

\sep
%

\begin{align*}%
{\displaystyle\int}
n  &  =\frac{1}{2}nn+\frac{1}{2}n.\\%
{\displaystyle\int}
nn  &  =\frac{1}{3}n^{3}+\frac{1}{2}nn+\frac{1}{6}n.\\%
{\displaystyle\int}
n^{3}  &  =\frac{1}{4}n^{4}+\frac{1}{2}n^{3}+\frac{1}{4}nn.\\%
{\displaystyle\int}
n^{4}  &  =\frac{1}{5}n^{5}+\frac{1}{2}n^{4}+\frac{1}{3}n^{3}\ast-\frac{1}%
{30}n.\\%
{\displaystyle\int}
n^{5}  &  =\frac{1}{6}n^{6}+\frac{1}{2}n^{5}+\frac{5}{12}n^{4}\ast-\frac
{1}{12}nn.\\%
{\displaystyle\int}
n^{6}  &  =\frac{1}{7}n^{7}+\frac{1}{2}n^{6}+\frac{1}{2}n^{5}\ast-\frac{1}%
{6}n^{3}\ast+\frac{1}{42}n.\\%
{\displaystyle\int}
n^{7}  &  =\frac{1}{8}n^{8}+\frac{1}{2}n^{7}+\frac{7}{12}n^{6}\ast-\frac
{7}{24}n^{4}\ast+\frac{1}{12}nn.\\%
{\displaystyle\int}
n^{8}  &  =\frac{1}{9}n^{9}+\frac{1}{2}n^{8}+\frac{2}{3}n^{7}\ast-\frac{7}%
{15}n^{5}\ast+\frac{2}{9}n^{3}\ast-\frac{1}{30}n.\\%
{\displaystyle\int}
n^{9}  &  =\frac{1}{10}n^{10}+\frac{1}{2}n^{9}+\frac{3}{4}n^{8}\ast-\frac
{7}{10}n^{6}\ast+\frac{1}{2}n^{4}\ast-\frac{3}{20}nn.\\%
{\displaystyle\int}
n^{10}  &  =\frac{1}{11}n^{11}+\frac{1}{2}n^{10}+\frac{5}{6}n^{9}\ast
-1n^{7}\ast+1n^{5}\ast-\frac{1}{2}n^{3}\ast+\frac{5}{66}n.
\end{align*}
\textsf{Indeed, a pattern can be seen in the progressions herein which can be
continued by means of this rule: Suppose that }$c$\textsf{\ is the value of
any power; then the sum of all }$n^{c}$\textsf{\ or }
\begin{align*}%
{\displaystyle\int}
n^{c}  &  =\frac{1}{c+1}n^{c+1}+\frac{1}{2}n^{c}+\frac{c}{2}An^{c-1}%
+\frac{c\cdot c-1\cdot c-2}{2\cdot3\cdot4}Bn^{c-3}\\
&  +\frac{c\cdot c-1\cdot c-2\cdot c-3\cdot c-4}{2\cdot3\cdot4\cdot5\cdot
6}Cn^{c-5}\\
&  +\frac{c\cdot c-1\cdot c-2\cdot c-3\cdot c-4\cdot c-5\cdot c-6}%
{2\cdot3\cdot4\cdot5\cdot6\cdot7\cdot8}Dn^{c-7}\ldots\text{\textsf{,}}%
\end{align*}
\textsf{where the value of the power }$n$\textsf{\ continues to decrease by
two until it reaches }$n$\textsf{\ or }$nn$\textsf{. The uppercase letters
}$A$\textsf{, }$B$\textsf{, }$C$\textsf{, }$D$\textsf{, etc., in order, denote
the coefficients of the final term of }$%
{\displaystyle\int}
nn,$\textsf{\ }$%
{\displaystyle\int}
n^{4},$\textsf{\ }$%
{\displaystyle\int}
n^{6},$\textsf{\ }$%
{\displaystyle\int}
n^{8}$\textsf{, etc., namely }
\[
A=\frac{1}{6},\emph{\ }B=-\frac{1}{30},C=\frac{1}{42},\ D=-\frac{1}%
{30}\text{.}%
\]
\textsf{These coefficients are such that, when arranged with the other
coefficients of the same order, they add up to unity: so, for }$D$\textsf{,
which we said signified }$-\frac{1}{30}$\textsf{, we have }
\[
\frac{1}{9}+\frac{1}{2}+\frac{2}{3}-\frac{7}{15}+\frac{2}{9}(+D)-\frac{1}%
{30}=1\text{.{}}%
\]
%

\sep

At this point we modern readers could conceivably exhibit great retrospective
prescience, anticipate Euler's broader context of $\sum_{i=1}^{n}f(i)$, for
which Bernoulli's claimed summation formula above provides test functions of
the form $f(x)=x^{c}$, and venture a rash generalization:
\[
\sum_{i=1}^{n}f(i)\approx C+\int^{n}f(x)dx+\frac{f(n)}{2}+A\frac{f^{\prime
}(n)}{2!}+B\frac{f^{\prime\prime\prime}(n)}{4!}+\cdots\text{\ }.
\]

This formula is what Euler discovered in the early 1730s (although he was
apparently unaware of Bernoulli's claim until later). Euler's summation
formula captures the delicate details of the general connection between
integration and discrete summation, and subsumes and resolves the two-thousand
year old quest for sums of powers formulas as a simple special case. In what
follows I will focus on just a few highlights from Euler.\smallskip

\section{The Basel problem and the summation formula}

``\emph{Euler calculated without any apparent effort, just as men breathe, as
eagles sustain themselves in the air.}'', Arago. \cite[p.\ 354]{young}

\smallskip

Around the year 1730, the 23-year old Euler, along with his frequent
correspondents Christian Goldbach (1690--1764) and Daniel Bernoulli
(1700--1782), developed ways to find increasingly accurate fractional or
decimal estimates for the sum of the reciprocal squares. But highly accurate
estimates were challenging, since the series converges very slowly. They were
likely trying to guess the exact value of the sum, hoping to recognize that
their approximations hinted something familiar, perhaps involving $\pi$, like
Leibniz's series, which had summed to $\pi/4$. Euler hit gold with the
discovery of his summation formula. One of its first major uses was in a
paper\footnote{E 47 in the Enestr\"{o}m Index \cite{archive}.} submitted to
the St.\ Petersburg Academy of Sciences on the 13th of October, 1735, in which
he approximated the sum correctly to twenty decimal places. Only seven and a
half weeks later Euler astonished his contemporaries with another
paper\footnote{E 41.}, solving the famous Basel problem by
demonstrating with a completely different method that the precise sum of the
series is $\pi^{2}/6$: \textquotedblleft\emph{Now, however, quite
unexpectedly, I have found an elegant formula for }$1+\frac{1}{4}+\frac{1}%
{9}+\frac{1}{16}+$\emph{\ etc., depending upon the quadrature of the circle
[i.e., upon }$\pi$\emph{]}\textquotedblright\ \cite[p.\ 261]{weil}. Johann
Bernoulli reacted \textquotedblleft\emph{And so is satisfied the burning
desire of my brother [Jakob] who, realizing that the investigation of the sum
was more difficult than anyone would have thought, openly confessed that all
his zeal had been mocked. If only my brother were alive now}\textquotedblright%
\ \cite[p.\ 345]{young}.

Much of Euler's \emph{Calculi Differentialis}, written two decades later, 
focused on the relationship between differential calculus and infinite series,
unifying his many discoveries in a single exposition. He devoted Chapters 5
and 6 of Part Two to the summation formula and a treasure trove of
applications. In Chapter 5 Euler derived his summation formula, analyzed the
generating function for Bernoulli numbers in terms of transcendental
functions, derived several properties of Bernoulli numbers, showed that they
grow supergeometrically, proved Bernoulli's formulas for sums of powers, and
found the exact sums of all infinite series of reciprocal even powers in terms
of Bernoulli numbers. Chapter 6 applied the summation formula to approximate
harmonic partial sums and the associated ``Euler'' constant $\gamma$, sums of
reciprocal powers, $\pi$, and sums of logarithms, leading to approximations
for large factorials and binomial coefficients.

I will guide the reader through a few key passages from the
translation. The reader may find more background, annotation, and exercises in
our book \cite{masters} or explore my more extensive translation on the web
\cite{instcalcdiff}. The passages below contain Euler's derivation,
the relation to Bernoulli numbers, application to reciprocal squares, and to
sums of logarithms, large factorials, and binomials, with mention of other
omitted passages. Each application uses the summation formula in a
fundamentally different way. The complete glory of Euler's chapters is
still available only in the original Latin \cite[volume 10]{euleroo} or an old
German translation \cite{eulertransl} (poorly printed in Fraktur); I 
encourage the reader to revel in the original.

\section{Euler's derivation}

Euler's derivation of his summation formula rests on two ideas. First, he used
Taylor series from calculus to relate the sum of the values of a function at 
finitely many
successive integers to similar sums involving the derivatives of the function.%

\sep

\begin{center}
\textsf{Leonhard Euler, from }

\textsf{Foundations of Differential Calculus}

\textsf{Part Two, Chapter 5}

\textsf{On Finding Sums of Series from the General Term}\smallskip
\end{center}

\textsf{105. Consider a series whose general term, belonging to the index }%
$x$\textsf{, is }$y$\textsf{, and whose preceding term, with index }%
$x-1$\textsf{, is }$v$\textsf{; because }$v$\textsf{\ arises from }%
$y$\textsf{, when }$x$\textsf{\ is replaced by }$x-1$\textsf{, one
has}\footnote{Euler expressed the value $v$ of his function at $x-1$ in terms
of its value $y$ at $x$ and the values of all its derivatives, also implicitly
evaluated at $x$. This uses Taylor series with increment $-1$. Of course 
he was tacitly assuming that this all makes sense, i.e., that his
function is infinitely differentiable, and that the Taylor series
converges and equals its intended value. Note also that the symbols $x$ and
$y$ are being used, respectively, to indicate the final value of an integer
index and the final value of the function evaluated there, as well as more
generally as a variable and a function of that variable. Today we would find
this much too confusing to dare write this way.}\textsf{\ }
\[
v=y-\frac{dy}{dx}+\frac{ddy}{2dx^{2}}-\frac{d^{3}y}{6dx^{3}}+\frac{d^{4}%
y}{24dx^{4}}-\frac{d^{5}y}{120dx^{5}}+\text{\textsf{\ etc.}}%
\]
\textsf{If }$y$\textsf{\ is the general term of the series }
\[
\left.  \setlength{\arraycolsep}{1pt}%
\begin{array}
[c]{ccccccccccccc}%
1 &  & 2 &  & 3 &  & 4 &  & \cdots &  & x-1 &  & x\\
a & + & b & + & c & + & d & + & \cdots & + & v & + & y
\end{array}
\right.
\]
\textsf{and if the term belonging to the index }$0$\textsf{\ is }$A$\textsf{,
then }$v$\textsf{, as a function of }$x$\textsf{, is the general term of the
series }
\[
\left.  \setlength{\arraycolsep}{1pt}%
\begin{array}
[b]{rrrrrrrrrrrrr}%
1 &  & 2 &  & 3 &  & 4 &  & 5 &  & \cdots &  & x\\
A & + & a & + & b & + & c & + & d & + & \cdots & + & v
\end{array}
\right.  \text{\textsf{,}}%
\]
\textsf{so if }$Sv$\textsf{\ denotes the sum of this series, then }%
$Sv=Sy-y+A$\textsf{.}

\textsf{106. Because }
\[
v=y-\frac{dy}{dx}+\frac{ddy}{2dx^{2}}-\frac{d^{3}y}{6dx^{3}}%
+\text{\textsf{\ etc.,}}%
\]
\textsf{one has, from the preceding, }
\[
Sv=Sy-S\frac{dy}{dx}+S\frac{ddy}{2dx^{2}}-S\frac{d^{3}y}{6dx^{3}}+S\frac
{d^{4}y}{24dx^{4}}-\text{\textsf{\ etc.,}}%
\]
\textsf{and, because }$Sv=Sy-y+A$\textsf{, }
\[
y-A=S\frac{dy}{dx}-S\frac{ddy}{2dx^{2}}+S\frac{d^{3}y}{6dx^{3}}-S\frac{d^{4}%
y}{24dx^{4}}+\text{\textsf{\ etc.,}}%
\]
\textsf{or equivalently }
\[
S\frac{dy}{dx}=y-A+S\frac{ddy}{2dx^{2}}-S\frac{d^{3}y}{6dx^{3}}+S\frac{d^{4}%
y}{24dx^{4}}-\text{\textsf{\ etc.}}%
\]
\textsf{Thus if one knows the sums of the series, whose general terms are
}$\frac{ddy}{dx^{2}},$\textsf{\ }$\frac{d^{3}y}{dx^{3}},$\textsf{\ }%
$\frac{d^{4}y}{dx^{4}},$\textsf{\ etc., one can obtain the summative term of
the series whose general term is }$\frac{dy}{dx}$\textsf{. The constant }%
$A$\textsf{\ must then be such that the summative term }$S\frac{dy}{dx}%
$\textsf{\ disappears when }$x=0$ ...%

\sep

Euler next applied this equation recursively, in \S 107--108,\ to demonstrate
how one can obtain individual sums of powers formulas, because in these cases
the derivatives will eventually vanish. He then continued with his second
idea, which produced the summation formula.%

\sep

\textsf{109. ... if one sets }$\frac{dy}{dx}=z$\textsf{, then }
\[
Sz=\int zdx+\frac{1}{2}S\frac{dz}{dx}-\frac{1}{6}S\frac{ddz}{dx^{2}}+\frac
{1}{24}S\frac{d^{3}z}{dx^{3}}-\text{\textsf{\ etc.,}}%
\]
\textsf{adding to it a constant value such that when }$x=0$\textsf{, the sum
}$Sz\ $\textsf{also vanishes. ... }\smallskip

\textsf{110. But if in the expressions above one substitutes the letter }%
$z$\textsf{\ in place of }$y$\textsf{, or if one differentiates the preceding
equation, which yields the same, one obtains }
\[
S\frac{dz}{dx}=z+\frac{1}{2}S\frac{ddz}{dx^{2}}-\frac{1}{6}S\frac{d^{3}%
z}{dx^{3}}+\frac{1}{24}S\frac{d^{4}z}{dx^{4}}-\text{\textsf{\ etc.;}}%
\]
\textsf{but using }$\frac{dz}{dx}$\textsf{\ in place of }$y$\textsf{\ one
obtains }
\[
S\frac{ddz}{dx^{2}}=\frac{dz}{dx}+\frac{1}{2}S\frac{d^{3}z}{dx^{3}}-\frac
{1}{6}S\frac{d^{4}z}{dx^{4}}+\frac{1}{24}S\frac{d^{5}z}{dx^{5}}%
-\text{\textsf{\ etc.}}%
\]
\textsf{... and so forth indefinitely.... }

\textsf{111. Now when these values for }$S\frac{dz}{dx}$\textsf{, }%
$S\frac{ddz}{dx^{2}}$\textsf{, }$S\frac{d^{3}z}{dx^{3}}$\textsf{\ are
successively substituted in the expression }
\[
Sz=\int zdx+\frac{1}{2}S\frac{dz}{dx}-\frac{1}{6}S\frac{ddz}{dx^{2}}+\frac
{1}{24}S\frac{d^{3}z}{dx^{3}}-\text{\textsf{\ etc.,}}%
\]
\textsf{one finds an expression for }$Sz$\textsf{, composed of the terms
}$\int zdx$\textsf{, }$z$\textsf{, }$\frac{dz}{dx}$\textsf{, }$\frac
{ddz}{dx^{2}}$\textsf{, }$\frac{d^{3}z}{dx^{3}}$\textsf{\ etc., whose
coefficients are easily obtained as follows. One sets }
\[
Sz=\int zdx+\alpha z+\frac{\beta dz}{dx}+\frac{\gamma ddz}{dx^{2}}%
+\frac{\delta d^{3}z}{dx^{3}}+\frac{\varepsilon d^{4}z}{dx^{4}}%
+\text{\textsf{\ etc.,}}%
\]
\textsf{and substitutes for these terms the values they have from the previous
series, yielding }
\[
\left.  \setlength{\arraycolsep}{1pt}%
\begin{array}
[c]{rrrrrrrrrrrrr}%
\int zdx & = & Sz & - & \frac{1}{2}S\frac{dz}{dx} & + & \frac{1}{6}S\frac
{ddz}{dx^{2}} & - & \frac{1}{24}S\frac{d^{3}z}{dx^{3}} & + & \frac{1}%
{120}S\frac{d^{4}z}{dx^{4}} & - & \text{\textsf{etc.}}\\
\vspace{-0.1in} &  &  &  &  &  &  &  &  &  &  &  & \\
\alpha z & = &  & + & \alpha S\frac{dz}{dx} & - & \frac{\alpha}{2}S\frac
{ddz}{dx^{2}} & + & \frac{\alpha}{6}S\frac{d^{3}z}{dx^{3}} & - & \frac{\alpha
}{24}S\frac{d^{4}z}{dz^{4}} & + & \text{\textsf{etc.}}\\
\vspace{-0.1in} &  &  &  &  &  &  &  &  &  &  &  & \\
\frac{\beta dz}{dx} & = &  &  &  &  & \beta S\frac{ddz}{dx^{2}} & - &
\frac{\beta}{2}S\frac{d^{3}z}{dx^{3}} & + & \frac{\beta}{6}S\frac{d^{4}%
z}{dx^{4}} & - & \text{\textsf{etc.}}\\
\vspace{-0.1in} &  &  &  &  &  &  &  &  &  &  &  & \\
\frac{\gamma ddz}{dx^{2}} & = &  &  &  &  &  &  & \gamma S\frac{d^{3}z}%
{dx^{3}} & - & \frac{\gamma}{2}S\frac{d^{4}z}{dx^{4}} & + &
\text{\textsf{etc.}}\\
\vspace{-0.1in} &  &  &  &  &  &  &  &  &  &  &  & \\
\frac{\delta d^{3}z}{dx^{3}} & = &  &  &  &  &  &  &  &  & \delta
\;S\frac{d^{4}z}{dx^{4}} & - & \text{\textsf{etc.}}\\
\vspace{-0.1in} &  &  &  &  &  &  &  &  &  &  &  & \\
&  &  &  &  &  & \text{\textsf{etc.}} &  &  &  &  &  &
\end{array}
\right.
\]
\textsf{Since these values, added together, must produce }$Sz$\textsf{, the
coefficients }$\alpha$\textsf{, }$\beta$\textsf{, }$\gamma$\textsf{, }$\delta
$\textsf{\ etc. are ... }

\textsf{112. ...}
\begin{gather*}
\alpha=\frac{1}{2}\text{,\thinspace}\beta=\frac{\alpha}{2}-\frac{1}{6}%
=\frac{1}{12}\text{,\thinspace}\gamma=\frac{\beta}{2}-\frac{\alpha}{6}%
+\frac{1}{24}=0,\\[0.05in]
\delta=\frac{\gamma}{2}-\frac{\beta}{6}+\frac{\alpha}{24}-\frac{1}{120}%
=-\frac{1}{720}\text{,\thinspace}\varepsilon=\frac{\delta}{2}-\frac{\gamma}%
{6}+\frac{\beta}{24}-\frac{\alpha}{120}+\frac{1}{720}=0\text{\textsf{\ etc.,}}%
\end{gather*}
\textsf{and if one continues in this fashion one finds that alternating terms
vanish.}

\section{Connection to Bernoulli numbers and sums of powers}

Before Euler showed how to apply his summation formula to derive new results, 
in \S 112--120\ he intensively studied the coefficients
$\alpha,\beta,\gamma,\ldots$, and discovered that their
generating function relates directly to the transcendental functions of
calculus, especially the cotangent. In particular, Euler proved that every 
second coefficient vanishes, and that those
that remain alternate in sign, by investigating a power series solution
to the  
differential equation satisfied by the cotangent function by dint of its 
derivative formula. Euler also explored 
number theoretic properties of the coefficients, including the growth and
prime factorizations of their numerators and denominators, some of which we
will see below.

\textbf{Caution}: In the process of distilling the summation formula in terms
of Bernoulli numbers, Euler switched the meaning of the Greek letters $\alpha
$, $\beta$, $\gamma$, $\delta$,..., and the formula now takes revised form:%

\sep

\textsf{121. ... If one finds the values of the [redefined] letters }$\alpha
$\textsf{, }$\beta$\textsf{, }$\gamma$\textsf{, }$\delta$\textsf{, etc.
according to this rule, which entails little difficulty in calculation, then
one can express the summative term of any series, whose general term }%
$=z$\textsf{\ corresponding to the index }$x$\textsf{, in the following
fashion: }
\begin{align*}
Sz  &  =\int zdx+\frac{1}{2}z+\frac{\alpha dz}{1\cdot2\cdot3dx}-\frac{\beta
d^{3}z}{1\cdot2\cdot3\cdot4\cdot5dx^{3}}+\frac{\gamma d^{5}z}{1\cdot
2\cdots7dx^{5}}\\
&  -\frac{\delta d^{7}z}{1\cdot2\cdots9dx^{7}}+\frac{\varepsilon d^{9}%
z}{1\cdot2\cdots11dx^{9}}-\frac{\zeta d^{11}z}{1\cdot2\cdots13dx^{11}%
}+\text{\textsf{etc. }}\ldots
\end{align*}

\textsf{122. These numbers have great use throughout the entire theory of
series. First, one can obtain from them the final terms in the sums of even
powers, for which we noted above (in \S 63 of part one) that one cannot obtain
them, as one can the other terms, from the sums of earlier powers. For the
even powers, the last terms of the sums are products of }$x$\textsf{\ and
certain numbers, namely for the 2nd, 4th, 6th, 8th, etc., }$\frac{1}{6}%
,$\textsf{\ }$\frac{1}{30},$\textsf{\ }$\frac{1}{42},$\textsf{\ }$\frac{1}%
{30}$\textsf{\ etc. with alternating signs. But these numbers arise from the
values of the letters }$\alpha$\textsf{, }$\beta$\textsf{, }$\gamma$\textsf{,
}$\delta$\textsf{, etc., which we found earlier, when one divides them by the
odd numbers }$3,$\textsf{\ }$5,$\textsf{\ }$7,$\textsf{\ }$9$\textsf{, etc.
These numbers are called the Bernoulli numbers after their discoverer Jakob
Bernoulli, and they are }
\[
\left.  \setlength{\arraycolsep}{1pt}%
\begin{array}
[c]{rccclrrcccccl}%
\frac{\alpha}{3} & = & \frac{1}{6} & = & \mathfrak{A} & \qquad\text{\quad} &
\frac{\iota}{19} & = & \frac{43867}{798} & = & \mathfrak{I} &  & \\
\vspace{-0.1in} &  &  &  &  &  &  &  &  &  &  &  & \\
\frac{\beta}{5} & = & \frac{1}{30} & = & \mathfrak{B} &  & \frac{\chi}{21} &
= & \frac{174611}{330} & = & \mathfrak{K} & = & \frac{283\cdot617}{330}\\
\vspace{-0.1in} &  &  &  &  &  &  &  &  &  &  &  & \\
\frac{\gamma}{7} & = & \frac{1}{42} & = & \mathfrak{C} &  & \frac{\lambda}{23}
& = & \frac{854513}{138} & = & \mathfrak{L} & = & \frac{11\cdot131\cdot
593}{2\cdot3\cdot23}\\
\vspace{-0.1in} &  &  &  &  &  &  &  &  &  &  &  & \\
\frac{\delta}{9} & = & \frac{1}{30} & = & \mathfrak{D} &  & \frac{\mu}{25} &
= & \frac{236364091}{2730} & = & \mathfrak{M} &  & \\
\vspace{-0.1in} &  &  &  &  &  &  &  &  &  &  &  & \\
\frac{\varepsilon}{11} & = & \frac{5}{66} & = & \mathfrak{E} &  & \frac{\nu
}{27} & = & \frac{8553103}{6} & = & \mathfrak{N} & = & \frac{13\cdot657931}%
{6}\\
\vspace{-0.1in} &  &  &  &  &  &  &  &  &  &  &  & \\
\frac{\zeta}{13} & = & \frac{691}{2730} & = & \mathfrak{F} &  & \frac{\xi}{29}
& = & \frac{23749461029}{870} & = & \mathfrak{O} &  & \\
\vspace{-0.1in} &  &  &  &  &  &  &  &  &  &  &  & \\
\frac{\eta}{15} & = & \frac{7}{6} & = & \mathfrak{G} &  & \frac{\pi}{31} & = &
\frac{8615841276005}{14322} & = & \mathfrak{P} &  & \\
\vspace{-0.1in} &  &  &  &  &  &  &  &  &  &  &  & \\
\frac{\theta}{17} & = & \frac{3617}{510} & = & \mathfrak{H} &  &  &  &
\text{\textsf{etc.}} &  &  &  &
\end{array}
\right.  \mathsf{\smallskip}%
\]
%

\sep

Euler's very first application of the Bernoulli numbers, in \S 124--125, was to
solve a problem dear to his heart, determining the precise sums of all
infinite series of reciprocal even powers. His result (using today's notation
$\mathfrak{A=}B_{2},$ $\mathfrak{B=}-B_{4},$ $\mathfrak{C=}B_{6},$ $\ldots$)
was:
\[
\sum_{i=1}^{\infty}\frac{1}{i^{2n}}=\frac{\left(  -1\right)  ^{n+1}%
B_{2n}2^{2n-1}}{\left(  2n\right)  !}\pi^{2n}\text{ for all }n\geq1.
\]
Because these sums approach one as $n$ grows, he also obtained, in \S 129, an
asymptotic understanding of how Bernoulli numbers grow:
\[
\frac{B_{2n+2}}{B_{2n}}\approx-\frac{\left(  2n+2\right)  \left(  2n+1\right)
}{4\pi^{2}}\approx-\frac{n^{2}}{\pi^{2}}\text{.}%
\]
Thus he commented that they ``\emph{form a highly diverging sequence, which
grows more strongly than any geometric sequence of growing terms}''.

This completed Euler's analysis of the Bernoulli numbers. Now he was
ready to turn his summation formula towards applications. He ended Chapter 5
with applications in which the summation formula is finite (\S 131), including
that of a pure power function, which proved\ the formulas for sums of
powers discovered by Bernoulli (\S 132).

\section{``Until it begins to diverge''}

Chapter 6 applies the summation formula to make approximations even when it
diverges, which it does in almost all interesting situations.%

\sep

\begin{center}
\textsf{Part Two, Chapter 6}

\textsf{On the summing of progressions via infinite series\smallskip}
\end{center}

\textsf{140. The general expression, that we found in the previous chapter for
the summative term of a series, whose general term corresponding to the index
}$x$\textsf{\ is }$z$\textsf{, namely }
\[
Sz=\int zdx+\frac{1}{2}z+\frac{\mathfrak{A}dz}{1\cdot2dx}-\frac{\mathfrak{B}%
d^{3}z}{1\cdot2\cdot3\cdot4dx^{3}}+\frac{\mathfrak{C}d^{5}z}{1\cdot
2\cdots6dx^{5}}-\text{\textsf{\ etc.,}}%
\]
\textsf{actually serves to determine the sums of series, whose general terms
are integral rational functions}\footnote{By this he means polynomials.}%
\textsf{\ of the index }$x$\textsf{, because in these cases one eventually
arrives at vanishing differentials. On the other hand, if }$z$\textsf{\ is not
such a function of }$x$\textsf{, then the differentials continue without end,
and there results an infinite series that expresses the sum of the given
series up to and including the term whose index }$=x$\textsf{. The sum of the
series, continuing without end, is thus given by taking }$x=\infty$\textsf{,
and one finds in this way another infinite series equal to the original. ...}

\textsf{142. Since when a constant value is added to the sum, so that it
vanishes when }$x=0$\textsf{, the true sum is then found when }$x$\textsf{\ is
any other number, then it is clear that the true sum must likewise be given,
whenever a constant value is added that produces the true sum in any
particular case. Thus suppose it is not obvious, when one sets }$x=0$\textsf{,
what value the sum assumes and thus what constant must be used; one can
substitute other values for }$x$\textsf{, and through addition of a constant
value obtain a complete expression for the sum. Much will become clear from
the following. }%

\sep

For a particular choice of antiderivative $\int zdx$, the constant of interest 
is today called the \textquotedblleft Euler-Maclaurin
constant\textquotedblright\ for the function $z$ and a chosen antiderivative
$\int zdx$.

There follow Euler's \S 142a--144, in which he made the first application of
his summation formula to an infinite series, the diverging harmonic series
$\sum_{i=1}^{\infty}1/i$. For this series, the Euler-Maclaurin
constant in his summation formula will be the limiting difference between
$\sum_{i=1}^{x}1/i$ and $\ln x$. Today we call this particular number the 
\textquotedblleft Euler-Mascheroni constant,\textquotedblright and denote it by
$\gamma$. It is arguably the third most important constant in
mathematics after $\pi$ and $e$. Euler showed how to extract from the summation
formula an approximation of $\gamma$ accurate to 15 places and then easily
obtained the sum of the first thousand terms of the diverging harmonic series
to 13 places (see \cite{instcalcdiff}). In fact it is clear from what he
wrote that one could use his approach to approximate $\gamma$ to whatever
accuracy desired, and then apply the summation formula to find the value of
arbitrarily large finite harmonic sums to that same accuracy. I will discuss
in a moment the paradox that he can obtain arbitrarily accurate approximations
for the Euler-Maclaurin constant of a function and a chosen antiderivative
from a diverging summation!

We continue on to see exactly how Euler applied the summation formula to that
old puzzle, the Basel problem.%

\sep

\textsf{148. After considering the harmonic series we wish to turn to
examining the series of reciprocals of the squares, letting }
\[
s=1+\frac{1}{4}+\frac{1}{9}+\frac{1}{16}+\cdots+\frac{1}{xx}\text{.}%
\]
\textsf{Since the general term of this series is }$z=\frac{1}{xx}$\textsf{,
then }$\int zdx=\frac{-1}{x}$\textsf{, the differentials of }$z$\textsf{\ are
}
\[
\frac{dz}{2dx}=-\frac{1}{x^{3}},\quad\frac{ddz}{2\cdot3dx^{2}}=\frac{1}{x^{4}%
}\text{,}\quad\frac{d^{3}z}{2\cdot3\cdot4dx^{3}}=-\frac{1}{x^{5}}%
\quad\text{\textsf{\ etc.,}}%
\]
\textsf{and the sum is }
\[
s=C-\frac{1}{x}+\frac{1}{2xx}-\frac{\mathfrak{A}}{x^{3}}+\frac{\mathfrak{B}%
}{x^{5}}-\frac{\mathfrak{C}}{x^{7}}+\frac{\mathfrak{D}}{x^{9}}-\frac
{\mathfrak{E}}{x^{11}}+\text{\textsf{\ etc.,}}%
\]
\textsf{where the added constant }$C$\textsf{\ is determined from one case in
which the sum is known. We therefore wish to set }$x=1$\textsf{. Since then
}$s=1$\textsf{, one has }
\[
C=1+1-\frac{1}{2}+\mathfrak{A}-\mathfrak{B+C-D+E-}\text{\ \textsf{etc.,}}%
\]
\textsf{but this series alone does not give the value of }$C$\textsf{, since
it diverges strongly.}%

\sep

On the face of it, these formulas seem both absurd and useless. The expression
Euler obtains for the Euler-Maclaurin constant $C$ is clearly a divergent
series. In fact the summation formula here diverges for every $x$ because of
the supergeometric growth established for Bernoulli numbers. Euler, however,
was not fazed: he has a plan for obtaining from such divergent series highly
accurate approximations for both very large finite and infinite series.

Euler's idea was to add up the terms in the summation formula only
``\emph{until it begins to diverge.}'' For those unfamiliar with the theory of
divergent series, this seems preposterous, but in fact it has sound
theoretical underpinnings. Euler's approach was
ultimately vindicated by the modern theory of asymptotic series
\cite{hairer,hardy,hildebra,knopp}. Euler himself was probably confident of
his results, despite the apparently shaky foundations in divergent series,
because he was continually checking and rechecking his answers by a variety of
theoretical and computational methods, boosting his confidence in their
correctness from many different angles. Let us see how Euler continues
analyzing the sum of reciprocal squares, begun above.

First he recalled that for this particular function, he already knew the value
of $C$ by other means.%

\sep

\textsf{Above we demonstrated that the sum of the series to infinity is
}$=\frac{\pi\pi}{6}$\textsf{, and therefore setting }$x=\infty$\textsf{, and
}$s=\frac{\pi\pi}{6}$\textsf{, we have }$C=\frac{\pi\pi}{6}$\textsf{, because
then all other terms vanish. Thus it follows that }
\[
1+1-\frac{1}{2}+\mathfrak{A-B+C-D+E-}\text{\emph{etc}\textsf{. }}=\frac{\pi
\pi}{6}\text{.}%
\]
%

\sep

Next Euler imagined that he didn't already know the sum of the infinite 
series of reciprocal squares, and approximated it using his summation formula, 
thereby performing a cross-check on both methods.%

\sep

\textsf{149. If the sum of this series were not known, then one would need to
determine the value of the constant }$C$\textsf{\ from another case, in which
the sum were actually found. To this aim we set }$x=10$\textsf{\ and actually
add up ten terms, obtaining}\footnote{Euler used commas (as still done in
Europe today) rather than points, for separating the integer and fractional
parts of a decimal.}\textsf{\ }
\[
\left.  \setlength{\arraycolsep}{1pt}%
\begin{tabular}
[c]{cccll}
& $s$ & $=$ & $1,549767731166540690$ & .\\
& \vspace{-0.1in} &  &  & \\
\textsf{Further, add} & $\frac{1}{x}$ & $=$ & $0,1$ & \\
& \vspace{-0.1in} &  &  & \\
\textsf{subtr.} & $\frac{1}{2xx}$ & $=$ & $0,005$ & \\
&  &  & $\overline{1,644767731166540690}$ & \\
& \vspace{-0.1in} &  &  & \\
\textsf{add} & $\frac{\mathfrak{A}}{x^{3}}$ & $=$ & $\underline
{0,000166666666666666}$ & \\
&  &  & $1,644934397833207356$ & \\
& \vspace{-0.1in} &  &  & \\
\textsf{subtr.} & $\frac{\mathfrak{B}}{x^{5}}$ & $=$ & $\underline
{0,000000333333333333}$ & \\
&  &  & $1,644934064499874023$ & \\
& \vspace{-0.1in} &  &  & \\
\textsf{add} & $\frac{\mathfrak{C}}{x^{7}}$ & $=$ & $\underline
{0,000000002380952381}$ & \\
&  &  & $1,644934066880826404$ & \\
& \vspace{-0.1in} &  &  & \\
\textsf{subtr.} & $\frac{\mathfrak{D}}{x^{9}}$ & $=$ &
\multicolumn{1}{c}{$\underline{0,000000000033333333}$} & \\
&  &  & $1,644934066847493071$ & \\
& \vspace{-0.1in} &  &  & \\
\textsf{add} & $\frac{\mathfrak{E}}{x^{11}}$ & $=$ & $\underline
{0,000000000000757575}$ & \\
&  &  & $1,644934066848250646$ & \\
& \vspace{-0.1in} &  &  & \\
\textsf{subtr.} & $\frac{\mathfrak{F}}{x^{13}}$ & $=$ & $\underline
{0,000000000000025311}$ & \\
&  &  & $1,644934066848225335$ & \\
& \vspace{-0.1in} &  &  & \\
\textsf{add} & $\frac{\mathfrak{G}}{x^{15}}$ & $=$ & $0,000000000000001166$ &
\\
& \vspace{-0.1in} &  &  & \\
\textsf{subtr.} & $\frac{\mathfrak{H}}{x^{17}}$ & $=$ &
\multicolumn{1}{r}{$71$} & \\
&  &  & $\overline{1,644934066848226430}$ & $=C.$%
\end{tabular}
\ \right.
\]
\textsf{This number is likewise the value of the expression }$\frac{\pi\pi}%
{6}$\textsf{, as one can find by calculation from the known value of }$\pi
$\textsf{. From this it is clear that, although the series }$\mathfrak{A}%
,$\textsf{\ }$\mathfrak{B},$\textsf{\ }$\mathfrak{C}$\textsf{, etc. diverges,
it nevertheless produces a true sum. }%

\sep

So, on the one hand the summation formula diverges for every $x$, and yet on
the other it can apparently be used to make very close approximations, in fact
arbitrarily close approximations, to $C$. How can this be?

Note that the terms Euler actually calculated appear to decrease rapidly,
giving the initial appearance, albeit illusory, that the series converges.
Examining the terms more closely, one can see evidence that their decrease is
slowing in a geometric sense, which hints at the fact that the series actually
diverges. Recall that Euler intended to sum only ``\emph{until it begins to
diverge.}'' How did he decide when this occurs? Notice that the series
alternates in sign, and thus the partial sums bounce back and forth,
at first apparently converging, then diverging as the terms themselves
eventually increase due to rapid growth of the Bernoulli numbers. Euler knew
to stop before the smallest bounce, with the expectation that the true sum he
sought lies between any partial sum and the next one, and is thus
bracketed most accurately if one stops just before the smallest term 
is included.

Much later, through the course of the nineteenth century, mathematicians would
wrestle with the validity, theory and usefulness of divergent series. Two
(divergent) views reflected this struggle, and exemplified the evolution of 
mathematics:

\begin{quotation}
``\emph{The divergent series are the invention of the devil, and it is a shame
to base on them any demonstration whatsoever. By using them, one may draw any
conclusion he pleases and that is why these series have produced so many
fallacies and so many paradoxes. ...I have become prodigiously attentive to
all this, for with the exception of the geometrical series, there does not
exist in all of mathematics a single infinite series the sum of which has been
determined rigorously. In other words, the things which are most important in
mathematics are also those which have the least foundation. ... That most of
these things are correct in spite of that is extraordinarily surprising. I am
trying to find a reason for this; it is an exceedingly interesting
question.}'', Niels Abel (1802--1829), 1826 \cite[p.\ 973f]{kline}%
.\smallskip\smallskip

\noindent``\emph{The series is divergent; therefore we may be able to do
something with it}'', Oliver Heaviside (1850--1925) \cite[p.\ 1096]{kline}.
\end{quotation}

Euler, long before this, was confident in proceeding according to his simple
dictum ``\emph{until it begins to diverge.}'' Indeed, it is astounding but
true that the summation formula does behave exactly as Euler used it for many
functions, including all the ones Euler was interested in. Today we know for
certain that such ``asymptotic series'' indeed bracket the desired answer,
and diverge more and more slowly for larger and larger values of $x$, making 
them extremely useful for 
approximations \cite{hardy,hildebra,knopp}\cite[chapter 47]{kline}.

One can explore the interplay of calculation versus accuracy achieved by
different choices for $x$. A smaller choice for $x$ will cause the summation
formula to begin to diverge sooner, and with a larger final bounce, yielding
less accuracy. On the other hand, 
a larger $x$ will ensure much more rapid achievement of a
given level of accuracy, and greater bounding accuracy (as small as desired)
for the answer, at the expense of having to compute a longer partial sum on
the left hand side to get the calculation off the ground. Asymptotic series
have become important in applications of differential equations to
physical problems \cite[chapter 47]{kline}.

Euler's next application of the summation formula, in \S 150--153, was to
approximate the sums of reciprocal odd powers. I remarked above that Euler's
very first application of the Bernoulli numbers was to determine the precise
sums of all infinite series of reciprocal even powers. 
Naturally he also would have loved to find formulas
for the reciprocal odd powers, and he explored this at length using the
summation formula. He produced highly accurate decimal approximations for sums
of reciprocal odd powers all the way through the fifteenth, hoping to see a
pattern analogous to the even powers, namely simple fractions times the
relevant power of $\pi$. \ The first such converging series is the sum of
reciprocal cubes $\sum_{i=1}^{\infty}1/i^{3}$. Euler computed it accurately to
seventeen decimal places. He was disappointed, however, to find that it is not
near an obvious rational multiple of $\pi^{3}$, nor did he have better luck
with the other odd powers. Even today we know little about these sums of odd
powers, although not for lack of trying.

Following this, in \S 154--156 Euler approximated $\pi$ to seventeen decimal
places using the inverse tangent and cotangent functions with the summation
formula. He actually expressed his own amazement that one can approximate
$\pi$ so accurately with such an easy calculation!

\section{How to determine (or not) factorials}

I will showcase next Euler's efficacious use of the summation formula to
approximate finite sums of logarithms, and thus by exponentiating, to
approximate very large factorials via the formula now known as Stirling's
asymptotic approximation. Notice particularly Euler's ingenious determination
of the Euler-Maclaurin constant in the summation formula, from Wallis'
infinite product for $\pi$.

I will also briefly explore whether the
summation formula can determine a factorial precisely, yielding surprising results.

To set the stage for Euler, notice that to estimate a factorial, one can
estimate $\log\left(  x!\right)  =\log1+\log2+\cdots+\log x$, using any base,
provided one also knows how to find antilogarithms.%

\sep

\textsf{157. Now we want to use for }$z$\textsf{\ transcendental functions of
}$x$\textsf{, and take }$z=lx$\textsf{\ for summing hyperbolic}\footnote{Euler
called \textquotedblleft hyperbolic\textquotedblright\ logarithm what we today
call \textquotedblleft natural\textquotedblright\ logarithm.}%
\textsf{\ logarithms, from which the ordinary can easily be recovered, so that
}
\[
s=l1+l2+l3+l4+\cdots+lx.
\]
\textsf{Because }$z=lx$\textsf{, }
\[
\int zdx=xlx-x\text{,}%
\]
\textsf{since its differential is }$dxlx$. \textsf{Then }
\begin{gather*}
\frac{dz}{dx}=\frac{1}{x},\ \frac{ddz}{dx^{2}}=-\frac{1}{x^{2}},\ \frac
{d^{3}z}{1\cdot2dx^{3}}=\frac{1}{x^{3}}\text{,}\\[0.05in]
\frac{d^{4}z}{1\cdot
2\cdot3dx^{4}}=-\frac{1}{x^{4}},\ \frac{d^{5}z}{1\cdot2\cdot3\cdot4dx^{5}%
}=\frac{1}{x^{5}},\ \text{\textsf{etc.}}%
\end{gather*}
\noindent\textsf{One concludes that }
\[
s=xlx-x+\frac{1}{2}lx+\frac{\mathfrak{A}}{1\cdot2x}-\frac{\mathfrak{B}}%
{3\cdot4x^{3}}+\frac{\mathfrak{C}}{5\cdot6x^{5}}-\frac{\mathfrak{D}}%
{7\cdot8x^{7}}+\text{\textsf{etc.}}+\text{\textsf{Const.}}%
\]
\textsf{But for this constant one finds, when one sets }$x=1$\textsf{, because
then }$s=l1=0$\textsf{, }
\[
C=1-\frac{\mathfrak{A}}{1\cdot2}+\frac{\mathfrak{B}}{3\cdot4}-\frac
{\mathfrak{C}}{5\cdot6}+\frac{\mathfrak{D}}{7\cdot8}-\text{\textsf{etc.,}}%
\]
\textsf{a series that, due to its great divergence, is quite unsuitable even
for determining the approximate value of }$C$\textsf{. \smallskip}

\textsf{158. Nevertheless we can not only approximate the correct value of
}$C$\textsf{, but can obtain it exactly, by considering Wallis's expression
for }$\pi$\textsf{\ provided in the }\emph{Introductio}\textsf{\ }\cite[volume
1, chapter 11]{introductio}\textsf{. This expression is }
\[
\frac{\pi}{2}=\frac{2\cdot2\cdot4\cdot4\cdot6\cdot6\cdot8\cdot8\cdot
10\cdot10\cdot12\cdot\text{\textsf{etc.}}}{1\cdot3\cdot3\cdot5\cdot
5\cdot7\cdot7\cdot9\cdot9\cdot11\cdot11\cdot\text{\textsf{etc.}}}%
\]
\textsf{Taking logarithms, one obtains from this }
\begin{gather*}
l\pi-l2=2l2+2l4+2l6+2l8+2l10+l12+\text{\textsf{etc.}}\\
-l1-2l3-2l5-2l7-2l9-2l11-\text{\textsf{etc.}}%
\end{gather*}
\textsf{Setting }$x=\infty$ \textsf{in the assumed series, we have }
\[
\setlength{\arraycolsep}{1pt}\renewcommand{\arraystretch}{1.5}%
\begin{array}
[c]{rcrcl}
&  & l1+l2+l3+l4+\cdots+lx & = & C+\left(  x+\frac{1}{2}\right)
lx-x\text{,}\\
\text{\textsf{thus}} & \  & l1+l2+l3+l4+\cdots+l2x & = & C+\left(  2x+\frac
{1}{2}\right)  l2x-2x\\
\text{\textsf{and}} &  & l2+l4+l6+l8+\cdots+l2x & = & C+\left(  x+\frac{1}%
{2}\right)  lx+xl2-x,\\
\text{\textsf{and therefore}} &  & l1+l3+l5+l7+\cdots+l\left(  2x-1\right)  &
= & xlx+\left(  x+\frac{1}{2}\right)  l2-x.
\end{array}
\]
\textsf{Thus because }
\[
\left.  \setlength{\arraycolsep}{1pt}%
\begin{array}
[c]{llllllllllll}%
l\frac{\pi}{2} & = & 2l2 & + & 2l4 & + & 2l6 & + & \cdots & + & 2l2x-l2x & \\
\vspace{-0.1in} &  &  &  &  &  &  &  &  &  &  & \\
& - & 2l1 & - & 2l3 & - & 2l5 & - & \cdots & - & 2l\left(  2x-1\right)  , &
\end{array}
\right.
\]
\textsf{letting }$x=\infty$ \textsf{yields }
\[
l\frac{\pi}{2}=2C+\left(  2x+1\right)  lx+2xl2-2x-l2-lx-2xlx-\left(
2x+1\right)  l2+2x,
\]
\textsf{and therefore }
\[
l\frac{\pi}{2}=2C-2l2,\text{\textsf{thus} }2C=l2\pi\text{\ \textsf{and }%
}C=\frac{1}{2}l2\pi\text{,}%
\]
\textsf{yielding the decimal fraction representation }
\[
C=0,9189385332046727417803297\text{,}%
\]
\textsf{thus simultaneously the sum of the series }
\[
1-\frac{\mathfrak{A}}{1\cdot2}+\frac{\mathfrak{B}}{3\cdot4}-\frac
{\mathfrak{C}}{5\cdot6}+\frac{\mathfrak{D}}{7\cdot8}-\frac{\mathfrak{E}%
}{9\cdot10}+\text{\textsf{etc.}}=\frac{1}{2}l2\pi\text{.}%
\]
\textsf{\smallskip}

\textsf{159. Since we now know the constant }$C=\frac{1}{2}l2\pi$\textsf{, one
can exhibit the sum of any number of logarithms from the series }$l1+l2+l3+$
\textsf{etc. If one sets }
\[
s=l1+l2+l3+l4+\cdots+lx\text{,}%
\]
\textsf{then }
\[
s=\frac{1}{2}l2\pi+\left(  x+\frac{1}{2}\right)  lx-x+\frac{\mathfrak{A}%
}{1\cdot2x}-\frac{\mathfrak{B}}{3\cdot4x^{3}}+\frac{\mathfrak{C}}{5\cdot
6x^{5}}-\frac{\mathfrak{D}}{7\cdot8x^{7}}+\text{\textsf{etc.}}%
\]
\textsf{if the proposed logarithms are hyperbolic; if however the proposed
logarithms are common, then one must take common logarithms also in the terms
}$\frac{1}{2}l2\pi+(x+\frac{1}{2})lx$ \textsf{for }$l2\pi$ \textsf{and }%
$lx$\textsf{, and multiply the remaining terms }
\[
-x+\frac{\mathfrak{A}}{1\cdot2x}-\frac{\mathfrak{B}}{3\cdot4x^{3}%
}+\text{\textsf{etc.}}%
\]
\textsf{of the series by }$0,434294481903251827=n$\textsf{. In this case the
common logarithms are }
\begin{align*}
l\pi &  =0,497149872694133854351268\\
l2  &  =\underline{0,301029995663981195213738}\\
l2\pi &  =0,798179868358115049565006\\
\frac{1}{2}l2\pi &  =0,399089934179057524782503.
\end{align*}

\begin{center}
\textsf{Example.}

\textsf{\emph{Find the sum of the first thousand common logarithms} }
\end{center}

\[
s=l1+l2+l3+\cdots+l1000.
\]

\textsf{So }$x=1000$\textsf{, and }
\[
\left.  \setlength{\arraycolsep}{1pt}
\begin{array}
[c]{rrrl}%
lx & = & 3,0000000000000 & \hspace*{-2pt},\\
\vspace{-0.1in} &  &  & \\
\text{\textsf{and thus }}xlx & = & 3000,0000000000000 & \\
\vspace{-0.1in} &  &  & \\
\frac{1}{2}lx & = & 1,5000000000000 & \\
\vspace{-0.1in} &  &  & \\
\frac{1}{2}l2\pi & = & 0,3990899341790 & \\
&  & \overline{3001,8990899341790} & \\
\vspace{-0.1in} &  &  & \\
\text{\textsf{subtr. }}nx & = & 434,2944819032518 & \\
&  & \overline{2567,6046080309272} & \hspace*{-2pt}.\\
\text{\textsf{Then \qquad}} &  &  & \\
\frac{n\mathfrak{A}}{1\cdot2x} & = & 0,0000361912068 & \\
\vspace{-0.1in} &  &  & \\
\text{\textsf{subtr. }}\frac{n\mathfrak{B}}{3\cdot4x^{3}} & = & \underline
{0,0000000000012} & \\
&  & 0,0000361912056 & \\
\vspace{-0.1in} &  &  & \\
\text{\textsf{add}} &  & \underline{2567,6046080309272} & \\
\text{\textsf{the sum sought} }s & = & 2567,6046442221328 & \hspace*{-2pt}.
\end{array}
\right.
\]

\textsf{Now because }$s$ \textsf{is the logarithm of a product of numbers }
\[
1\cdot2\cdot3\cdot4\cdot5\cdot6\cdots1000,
\]
\textsf{it is clear that this product, if one actually multiplies it out,
consists of }$2568$\textsf{\ figures, beginning with the figures }%
$4023872$\textsf{, with }$2561$\textsf{\ subsequent figures. }%

\sep

One wonders how accurate such factorial approximations from the summation
formula can actually be. Exponentiating Euler's summation formula above for a
sum of logarithms produces the Stirling asymptotic approximation:
\[
x!\approx\frac{\sqrt{2\pi x}x^{x}}{e^{x}}e^{\left(  \mathfrak{A/}\left(
1\cdot2x\right)  -\mathfrak{B/}\left(  3\cdot4x^{3}\right)  +\mathfrak{C/}%
\left(  5\cdot6x^{5}\right)  -\mathfrak{D/}\left(  7\cdot8x^{7}\right)
+\cdots\right)  }.
\]
Because the summation formula diverges for each $x$, the accuracy of this
approximation is theoretically limited. Yet the value sought always lies
between those of successive partial sums. Moreover, from the asymptotic growth
rate of Bernoulli numbers obtained earlier, we see that approximately the
first $\pi x$ terms in the exponent might be expected to decrease (recall that
for $x=1000$ Euler used only two terms), with divergence occurring after that.

To explore the accuracy achievable with this formula, let us denote by
$S(x,m)$ the approximation to $x!$ using the first $m$ terms in the exponent.
Note that the discussion above tells us to expect this approximation to
\textquotedblleft start to diverge\textquotedblright\ after using around
$x\pi$ terms in the exponent, i.e., near $S\left(  x,\pi x\right)  $.
Beginning modestly with $x=10$, calculations with Maple show that $3628800$ is
the only integer between $S(10,2)$ and $S(10,3)$, thus determining $10!$ on
the nose. So although the summation formula has limited accuracy, it suffices
to determine easily the integer $10!$ uniquely.

And for $x=50$, one finds that \newline%
$30414093201713378043612608166064768844377641568960512000000000000$ \newline
is the only integer between $S(50,26)$ and $S(50,27)$, thus producing all $65$
digits of $50!$. This striking accuracy, and using so few of the roughly $\pi
x$ terms that in each case we expect to provide ever better approximations,
leads us to ask:

\textbf{Question}: Can one obtain the exact value of any factorial this way?

There is an interplay here as $x$ grows. Certainly the exponent becomes more
accurately known for larger $x$, using a given number of terms, and moreover
even more precisely known from the diverging series generally, which improves
for around $x\pi$ terms. On the other hand, it is then being exponentiated,
and finally, multiplied times something growing, to produce the factorial
approximation. So it is not so clear whether the factorial itself will always
be sufficiently trapped to determine its integer value.

Continuing experimentally, let us compare with $x=100$. First, note that
$100!$ is approximately $9.\allowbreak33\times10^{157}$. Using the same number
of terms, $27$, as was needed above to determine uniquely all $65$ digits of
$50!$, one finds that $S(100,27)$ agrees with $100!$ for the first $82$
digits. Thus it is giving more digits than when $x=50$, but does not yet
determine all the digits of $100!$. Further calculation shows that $100!$ is
however the unique integer first bracketed by $S(100,74)$ and $S(100,75)$. In
fact, from above one expects improvement for $100\pi$ terms. While one still
seems to have lots of terms to spare, one worries that, as $x$ increases, with
the number of decreasing terms in the summation only increasing linearly with
$x$, i.e., as $\pi x$, the number of terms needed to bracket the factorial
uniquely may be growing faster than this. In particular, when one doubled $x$
from $50$ to $100$, the number of terms needed to determine the factorial
increased from $27$ to $75$, more than doubling.

Both my theoretical analysis and further Maple computations ultimately confirm
this fear, eventually answering the question in the negative. But the size of
the factorials that are actually uniquely determined as integers by Euler's
summation formula, before it finally cannot keep up with all the digits, is
staggering. For instance, Euler showed above that $1000!$ possesses $2568$ 
digits, of which he calculated the first seven. My theoretical analysis shows 
that the Stirling approximation based on Euler's summation formula 
will determine every one of those $2568$ digits before it diverges.

\section{Large binomials}

In our final excerpt, Euler applied the summation formula to estimate the size
of large binomial coefficients. I translate just one of his methods here, in
which he merged two summation series term by term. As a sample application,
Euler studied the ratio $\binom{100}{50}/2^{100}$, despite the huge size
of its parts, thus closely approximating the probability that if one tosses
$100$ coins, exactly equal numbers will land heads and tails.%

\sep

\textsf{160. By means of this summation of logarithms, one can approximate the
product of any number of factors, that progress in the order of the natural
numbers. This can be especially helpful for the problem of finding the middle
or largest coefficient of any power in the binomial }$(a+b)^{m}$\textsf{,
where one notes that, when }$m$ \textsf{is an odd number, one always has two
equal middle coefficients, which taken together produce the middle coefficient
of the next even power. Thus since the largest coefficient of any even power
is twice as large as the middle coefficient of the immediately preceding odd
power, it suffices to determine the middle largest coefficient of an even
power. Thus we have }$m=2n$\textsf{\ with middle coefficient expressed as }
\[
\frac{2n\left(  2n-1\right)  \left(  2n-2\right)  \left(  2n-3\right)
\cdots\left(  n+1\right)  }{1\cdot2\cdot3\cdot4\cdots n}.
\]
\textsf{Setting this }$=u$\textsf{, one has }
\[
u=\frac{1\cdot2\cdot3\cdot4\cdot5\cdots2n}{\left(  1\cdot2\cdot3\cdot4\cdots
n\right)  ^{2}}\text{,}%
\]
\textsf{and taking logarithms }
\begin{align*}
lu  &  =l1+l2+l3+l4+l5+\cdots\;l2n\\
&  -2l1-2l2-2l3-2l4-2l5-\cdots-2ln\text{.}%
\end{align*}
\textsf{\smallskip}

\textsf{161. The sum of hyperbolic logarithms is }
\begin{align*}
l1+l2+l3+l4+\cdots+l2n  &  =\frac{1}{2}l2\pi+\left(  2n+\frac{1}{2}\right)
ln+\left(  2n+\frac{1}{2}\right)  l2-2n\\
&  +\frac{\mathfrak{A}}{1\cdot2\cdot2n}-\frac{\mathfrak{B}}{3\cdot4\cdot
2^{3}n^{3}}+\frac{\mathfrak{C}}{5\cdot6\cdot2^{5}n^{5}}-\text{\textsf{etc.}}%
\end{align*}
\textsf{and }
\begin{gather*}
2l1+2l2+2l3+2l4+\cdots+2ln\\
=l2\pi+\left(  2n+1\right)  ln-2n+\frac{2\mathfrak{A}}{1\cdot2n}%
-\frac{2\mathfrak{B}}{3\cdot4n^{3}}+\frac{2\mathfrak{C}}{5\cdot6n^{5}%
}-\text{\textsf{etc.}}%
\end{gather*}
\textsf{Subtracting this expression from the former yields }
\begin{gather*}
lu=-\frac{1}{2}l\pi-\frac{1}{2}ln+2nl2+\frac{\mathfrak{A}}{1\cdot2\cdot
2n}-\frac{\mathfrak{B}}{3\cdot4\cdot2^{3}n^{3}}+\frac{\mathfrak{C}}%
{5\cdot6\cdot2^{5}n^{5}}-\text{\textsf{etc.}}\\
-\frac{2\mathfrak{A}}{1\cdot2n}+\frac{2\mathfrak{B}}{3\cdot4n^{3}}%
-\frac{2\mathfrak{C}}{5\cdot6n^{5}}+\text{\textsf{etc.,}}%
\end{gather*}
\textsf{and collecting terms in pairs\footnote{Note that Euler's notation
leaves us to keep track of the scope of the square root symbol.} }
\[
lu=l\frac{2^{2n}}{\surd n\pi}-\frac{3\mathfrak{A}}{1\cdot2\cdot2n}%
+\frac{15\mathfrak{B}}{3\cdot4\cdot2^{3}n^{3}}-\frac{63\mathfrak{C}}%
{5\cdot6\cdot2^{5}n^{5}}+\frac{255\mathfrak{D}}{7\cdot8\cdot2^{7}n^{7}%
}-\text{\textsf{etc.}}%
\]

\textsf{...}

\textsf{162. ...}

\begin{center}
\textsf{\ Second Example }

\textsf{\emph{Find the ratio of the middle term of the binomial} }%
$(1+1)^{100}$ \textsf{\emph{to the sum} }$2^{100}$ \textsf{\emph{of all the
terms.}\ }
\end{center}

\textsf{For this we wish to use the formula we found first, }
\[
lu=l\frac{2^{2n}}{\surd n\pi}-\frac{3\mathfrak{A}}{1\cdot2\cdot2n}%
+\frac{15\mathfrak{B}}{3\cdot4\cdot2^{3}n^{3}}-\frac{63\mathfrak{C}}%
{5\cdot6\cdot2^{5}n^{5}}+\text{\textsf{etc.,}}%
\]
\textsf{from which, setting }$2n=m$\textsf{, in order to obtain the power
}$(1+1)^{m}$\textsf{, and after substituting the values of the letters
}$\mathfrak{A},$ $\mathfrak{B},$ $\mathfrak{C},$ $\mathfrak{D}$ \textsf{etc.,
one has }
\[
lu=l\frac{2^{m}}{\surd\frac{1}{2}m\pi}-\frac{1}{4m}+\frac{1}{24m^{3}}-\frac
{1}{20m^{5}}+\frac{17}{112m^{7}}-\frac{31}{36m^{9}}+\frac{691}{88m^{11}%
}-\text{\textsf{etc.}}%
\]
\textsf{Since the logarithms here are hyperbolic, one multiplies by }
\[
k=0,434294481903251\text{,}%
\]
\textsf{in order to change to tables, yielding }
\[
lu=l\frac{2^{m}}{\surd\frac{1}{2}m\pi}-\frac{k}{4m}+\frac{k}{24m^{3}}-\frac
{k}{20m^{5}}+\frac{17k}{112m^{7}}-\frac{31k}{36m^{9}}+\text{\textsf{etc.,}}%
\]
\textsf{Now since }$u$ \textsf{is the middle coefficient, the ratio sought is
}$2^{m}:u$\textsf{, and }
\[
l\frac{2^{m}}{u}=l\surd\frac{1}{2}m\pi+\frac{k}{4m}-\frac{k}{24m^{3}}+\frac
{k}{20m^{5}}-\frac{17k}{112m^{7}}+\frac{31k}{36m^{9}}-\frac{691k}{88m^{11}%
}+\text{\textsf{etc.}}%
\]
\textsf{Now, since the exponent }$m=100$\textsf{, }
\[
\frac{k}{m}=0,0043429448\text{,\quad}\frac{k}{m^{3}}=0,0000004343,\quad
\frac{k}{m^{5}}=0,0000000000,
\]
\textsf{yielding }
\[
\left.  \setlength{\arraycolsep}{1.0pt}
\begin{array}
[c]{rrl}%
\frac{k}{4m}= & 0,0010857362 & \\
\vspace{-0.1in} &  & \\
\frac{k}{24m^{3}}= & \underline{0,0000000181} & \\
\vspace{-0.1in} &  & \\
& \underline{0,0010857181} & .\\
\vspace{-0.1in} &  & \\
\text{\textsf{Further }}l\pi= & 0,4971498726 & \\
\vspace{-0.1in} &  & \\
l\frac{1}{2}m= & \underline{1,6989700043} & \\
\vspace{-0.1in} &  & \\
l\frac{1}{2}m\pi= & \underline{2,1961198769} & \\
\vspace{-0.1in} &  & \\
l\surd\frac{1}{2}m\pi= & 1,0980599384 & \\
\vspace{-0.1in} &  & \\
\frac{k}{4m}-\frac{k}{24m^{3}}+\text{\textsf{etc.}}= & \underline
{0,0010857181} & \\
\vspace{-0.1in} &  & \\
& 1,0991456565 & =l\frac{2^{100}}{u}.
\end{array}
\right.
\]
\textsf{Thus }$\frac{2^{100}}{u}=12,56451$\textsf{, and the middle term in the
expanded power }$(1+1)^{m}$ \textsf{\ is to the sum of all the terms }%
$2^{100}$ \textsf{as }$1$\textsf{\ is to }$12,56451$\textsf{. }%

\sep

So the probability that $100$ coin tosses will result in exactly $50$ each 
of heads and tails is between one in twelve and one in thirteen.

\end{document}